\documentclass[a4paper,11pt]{article}
\usepackage{amsmath,amssymb,color}
\usepackage{graphicx}
\usepackage{amsthm,bm,comment}
\usepackage{comment}
\usepackage{color}
\usepackage{natbib}
\usepackage{tikz} 
\usetikzlibrary{intersections, calc} 
\usepackage{setspace}
\definecolor{red}{rgb}{1,0,0} 
\definecolor{green}{rgb}{0,1,0} 
\definecolor{blue}{rgb}{0,0,1} 
\definecolor{purple}{rgb}{0.5,0,1} 

\theoremstyle{plain}
\newtheorem{defi}{Definition}
\newtheorem{thm}{Theorem}

\title{
	Copula bounds for circular data}
\author{
	Hiroaki Ogata\footnote{E-mail: hiroakiogata@tmu.ac.jp}
	\\
	{\scriptsize Faculty of Economics and Business Administration, Tokyo Metropolitan University} \\
	{\scriptsize 1-1 Minami-Osawa, Hachioji-shi, Tokyo, Japan 192-0397} 
}
\date{\today}
\begin{document}
	\maketitle 
	
\begin{abstract}
We propose the extension of Fr\'{e}chet-Hoeffding copula bounds for circular data. The copula is a powerful tool for describing the dependency of random variables. In two dimensions, the Fr\'{e}chet-Hoeffding upper (lower) bound indicates the perfect positive (negative) dependence between two random variables. However, for circular random variables, the usual concept of dependency is not accepted because of their periodicity. In this work, we redefine Fr\'{e}chet-Hoeffding bounds and consider modified Fr\'{e}chet and Mardia families of copulas for modelling the dependency of two circular random variables. Simulation studies are also given to demonstrate the behavior of the model.  \\  
{\bf keywords} : Circular statistics ; copulas ; Fr\'{e}chet-Hoeffding bounds
\end{abstract}

	\section{Introduction}\label{sec:Introduction}
	Circular statistics deals with observations that are represented as points on a unit circle. 
	Directional data is a typical example of this, because a direction is represented as an angle from a certain zero direction. Common examples often cited are: wind direction, river flow direction, and the direction of migrating birds in flight. For a comprehensive explanation of circular statistics, we refer to \citet{F95}, \citet{MJ99}, \citet{JS01}, to name a few. 
	
	If the observation is a pair of circular data, usually denoted by $ (\theta, \phi) \in \mathbf{C}^2=[0,2\pi)^2 $, it is represented as a point on a torus. When we have bivariate data, we can investigate the relationship between them. Although the Pearson correlation coefficient is the representative measure of a relationship between two random variables, it does not work for circular data because of its periodic structure. Many measures of association for circular data have been proposed so far. Section 11.2.2. of \citet{MJ99} provides a summary of circular-circular correlations. 
	
	Another powerful tool for describing the dependency between two random variables is a copula. It is a function rather than a single value and provides much more information regarding dependency. Mathematically speaking, the copula is a joint distribution function on $\mathbf{I}^2=[0,1]^2 $ whose margins are uniform. For a circular version, \citet{JPK15} introduced a circular analogue of copula density, called ``circulas,'' assuming the existence of density function. In this paper, we consider a copula function---not a copula density---for circular data. Because we do not assume the existence of the density function, we can deal with singular copulas. We give special consideration to circular analogues of Fr\'{e}chet-Hoeffding copula bounds, which are examples of singular copulas. Due to the arbitrary nature of the origin in the circular variable, they are not uniquely specified.   
	
	The remainder of the paper is organized as follows: Section \ref{sec:Equivalence} provides the definition of the equivalence class of circular copula functions from the aspect of the arbitrariness of the origin point. Section \ref{sec:bound} introduces circular analogues of Fr\'{e}chet-Hoeffding copula bounds. We prove that the following (i) and (ii) are equivalent. 
	\begin{enumerate}
		\renewcommand{\labelenumi}{(\roman{enumi})}
		\setlength{\itemsep}{0cm}
		\item\label{enu:upper bound} The bivariate circular random variables have the circular Fr\'{e}chet-Hoeffding upper (lower) copula bounds.
		\item\label{enu:nonincreasing} Support of the bivariate circular random variables is nondecreasing  (nonincreasing) in the sense of $ \mod 2\pi $.
	\end{enumerate}
	Section \ref{sec:Monte Carlo Simulations and Numerical Examples} gives a Monte Carlo simulation. We generate observations from the circular version of the Mardia copula family, and investigate their behavior. Section \ref{sec:Summary and Conclusions} provides our summary and conclusion. 
	
	\section{Equivalence Class of Circular Copula Functions}\label{sec:Equivalence}
	A copula is the joint distribution function on $ \mathbf{I}^2=[0,1]^2 $ whose margins are uniform. Sklar's theorem insists any joint distribution function $ H(x,y) $ is written by the copula function $ C(u,v) $ and their marginal distribution function $ F(x) $, $ G(y) $, like
	\begin{align*}
	H(x,y)=C(F(x),G(y)). 
	\end{align*}
	In the case of circular random variables, marginal and joint distribution functions depend on the choice of the zero direction. However, we should not consider the difference caused solely by the difference of the zero direction. This section gives the concept of the equivalence class of circular copula functions in the sense of the arbitrariness of the choice of the zero direction.  
	
	First, let us give the definitions of circular distribution functions, quasi-inverses of circular distribution functions and circular joint distribution functions following the style of Definitions 2.3.1, 2.3.6, and 2.3.2. in \citet{N06}.
	\begin{defi}
		A \textit{circular distribution function} is a function $ F $ with domain $ \mathbf{C}=[0, 2\pi) $ such that
		\begin{enumerate}
			\item $ F $ is nondecreasing,
			\item $ F(0) = 0 $ and $ \lim_{\theta \nearrow 2\pi} F(\theta) =1 $. 
		\end{enumerate}
	\end{defi}
	\begin{defi}
		Let $ F $ be a circular distribution function. Then a \textit{quasi-inverse} of $ F $ is any function $ F^{(-1)} $ with domain $ \mathbf{I} $ such that 
		\begin{align*}
		F^{(-1)}(u)=\inf \{ \theta | F(\theta)\ge u\}, 
		\end{align*}
	\end{defi}
	\begin{defi} 
		A \textit{circular joint distribution function} is a function $ H $ with domain $ \mathbf{C}^2 $ such that
		\begin{enumerate}
			\item $ H(\theta_2,\phi_2)-H(\theta_2,\phi_1)-H(\theta_1,\phi_2)+H(\theta_1,\phi_1) \ge 0 $ for any rectangle $ [\theta_1,\theta_2]\times[\phi_1,\phi_2]  \in \mathbf{C}^2$,
			\item $  H(\theta,0) =  H(0,\phi) = 0 $ and $ \lim_{\theta \nearrow 2\pi}\lim_{\phi \nearrow 2\pi} H(\theta,\phi)=1 $. 
		\end{enumerate}
	\end{defi}
	Now let us extend the domain of $ F $ to $ \tilde{\mathbf{C}}=[0,4\pi) $ in the following way
	\begin{align*}
	\tilde{F}(\theta)=\left\{
	\begin{array}{ll}
	F(\theta) & (0 \le \theta < 2\pi) \\
	F(\theta-2\pi)+1 & (2\pi \le \theta < 4\pi) 
	\end{array}
	\right.
	\end{align*}
	If we change the zero direction to $ \alpha \in \mathbf{C} $, its distribution function becomes
	\begin{align*}
	F_{\alpha}(\theta)=\tilde{F}(\theta+\alpha)-\tilde{F}(\alpha). 
	\end{align*}
	The choice of $ \alpha $ should be arbitrary, so we can define the equivalence class of circular distribution functions $ \{F_\alpha | \alpha \in \mathbf{C} \} $. 
	
	Similarly, we extend the domain of $ H $ to $ \tilde{\mathbf{C}}^2 $ in the following way
	\begin{align*}
	&\tilde{H}(\theta,\phi) \\
	=&\left\{
	\begin{array}{ll}
	H(\theta,\phi) & (0 \le \theta < 2\pi, 0 \le \phi < 2\pi) \\
	F(\theta)+H(\theta,\phi-2\pi) & (0 \le \theta < 2\pi, 2\pi \le \phi < 4\pi) \\
	G(\phi)+H(\theta-2\pi,\phi) & (2\pi \le \theta < 4\pi, 0 \le \phi < 2\pi) \\
	1+F(\theta-2\pi)+G(\phi-2\pi)+H(\theta-2\pi,\phi-2\pi) & (2\pi \le \theta < 4\pi, 2\pi \le \phi < 4\pi) 
	\end{array}
	\right.
	\end{align*}
	If we change the zero directions to $ (\alpha,\beta) \in \mathbf{C}^2 $, its joint distribution function becomes
	\begin{align*}
	H_{\alpha,\beta}(\theta,\phi)=\tilde{H}(\theta+\alpha,\phi+\beta)-\tilde{H}(\alpha,\phi+\beta)-\tilde{H}(\theta+\alpha,\beta)+\tilde{H}(\alpha,\beta). 
	\end{align*}
	The choice of $ (\alpha,\beta) $ should be arbitrary, so we can define the equivalence class of circular joint distribution functions $ \{H_{\alpha,\beta} | (\alpha,\beta) \in \mathbf{C}^2 \} $. 
	
	From Sklar's theorem, the circular copula function is defined by
	\begin{align}\label{eq:copula}
	C_{\alpha,\beta}(u,v)=H_{\alpha,\beta}(F_{\alpha}^{(-1)}(u),G_{\beta}^{(-1)}(v)).
	\end{align}
	The choice of zero directions $ (\alpha,\beta) $ should be arbitrary, so we can define the equivalence class of circular copula functions $ \{ C_{\alpha,\beta} | (\alpha,\beta) \in \mathbf{C}^2 \} $. 

	\section{Circular Fr\'{e}chet-Hoeffding copula bounds}\label{sec:bound}
	A copula is known to have both upper and lower bounds. That is, for every $ (u,v) \in \mathbf{I}^2 $, 
	\begin{align*}
	W(u,v):=\max(u+v-1,0) \le C(u,v) \le \min(u,v)=:M(u,v). 
	\end{align*} 
	The bounds $ M(u,v) $ and $ W(u,v) $ are themselves copulas and are called \textit{Fr\'{e}chet-Hoeffding upper bounds} and \textit{Fr\'{e}chet-Hoeffding lower bounds}, respectively. 
	
\begin{figure}[h] 
		\centering
			\includegraphics[width=60mm,height=60mm]{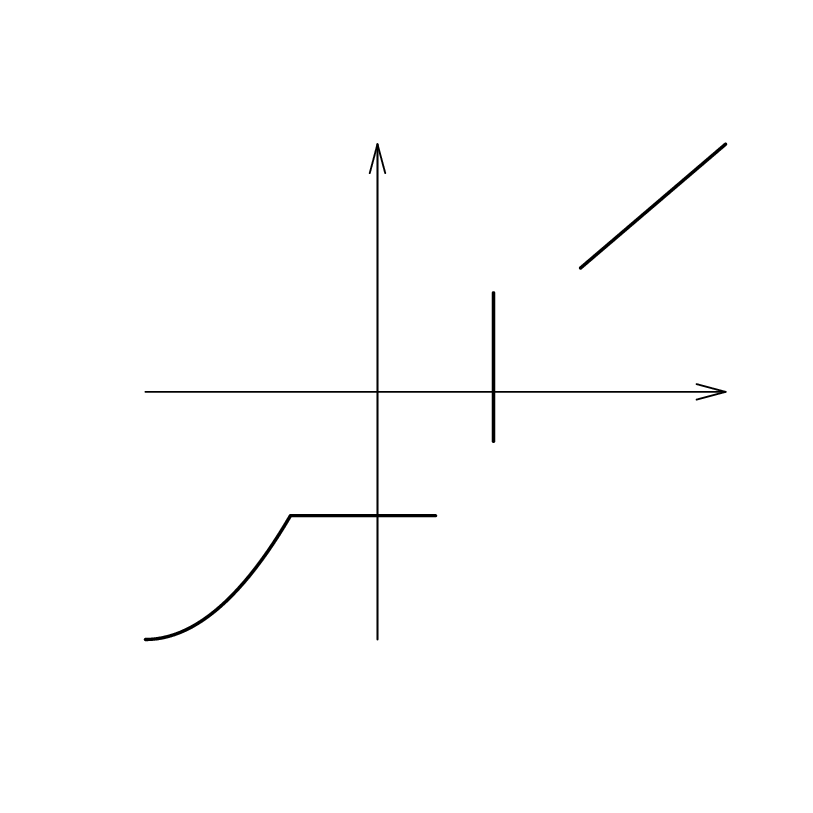}
		\caption{Nondecreasing set}
		\label{fig:nondecreasing}
	\end{figure}

	Here, we introduce the concept of a nondecreasing (nonincreasing) set in $ \bar{\mathbf{R}}^2=(-\infty,\infty)^2 $. 
	\begin{defi}[Definition 2.5.1 in \citet{N06}]
		A subset $ S $ of $ \bar{\mathbf{R}}^2 $ is \textit{nondecreasing} if for any $ (x,y) $ and $ (u,v) $ in $ S $, $ x < u $ implies $ y \le v $. Similarly, a subset $ S $ of $ \bar{\mathbf{R}}^2 $ is \textit{nonincreasing} if for any $ (x,y) $ and $ (u,v) $ is $ S $, $ x < u $ implies $ y \ge v $. 
	\end{defi}
	The example of a nondecreasing set is given in Figure \ref{fig:nondecreasing}.
	If the pair of random variables $ (X,Y) $ have a Fr\'{e}chet-Hoeffding upper (lower) bound, then the support of $ (X,Y) $ is nondecreasing (nonincreasing). See Theorems 2.5.4 and 2.5.5 of \citet{N06}. In this sense, a Fr\'{e}chet-Hoeffding upper (lower) bound indicates the perfect positive (negative) dependence between $ X $ and $ Y $. 
	
\begin{figure}[h] 
		\centering
				\includegraphics[width=50mm,height=50mm]{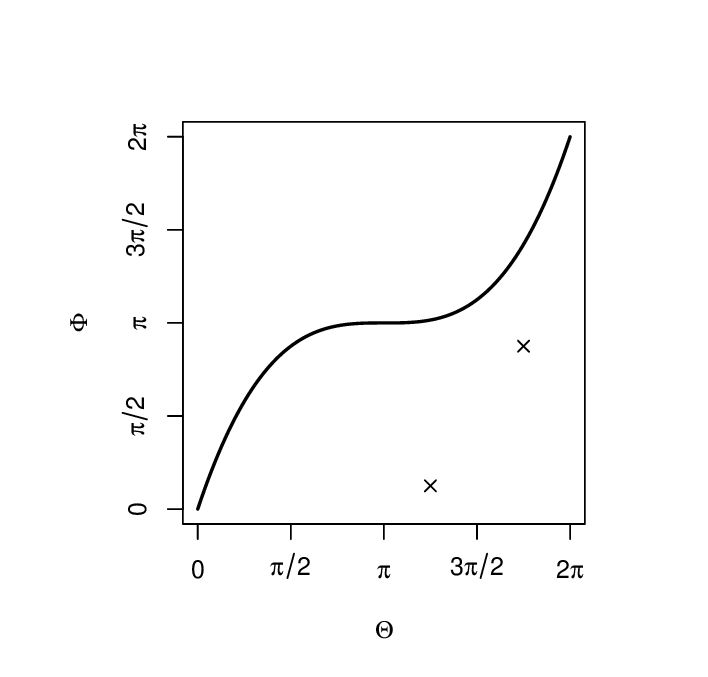}
				\includegraphics[width=50mm,height=50mm]{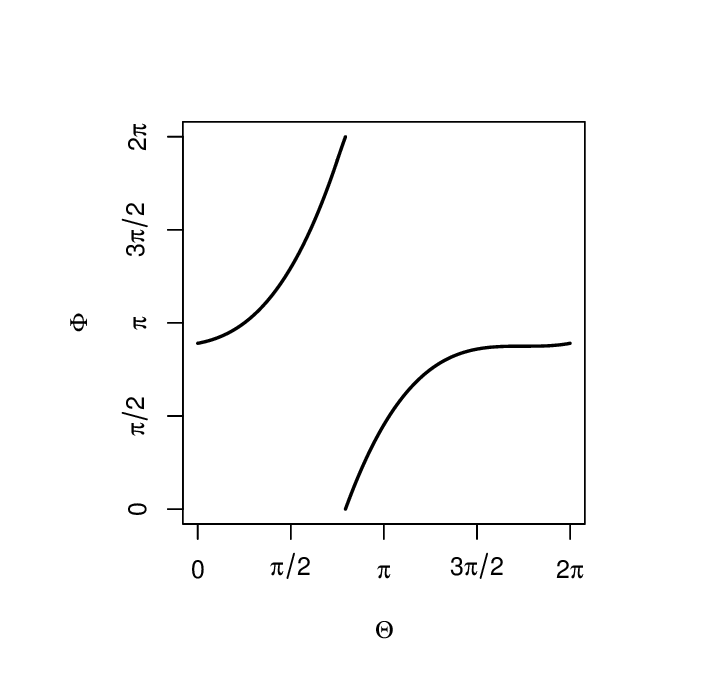}
				\includegraphics[width=50mm,height=50mm]{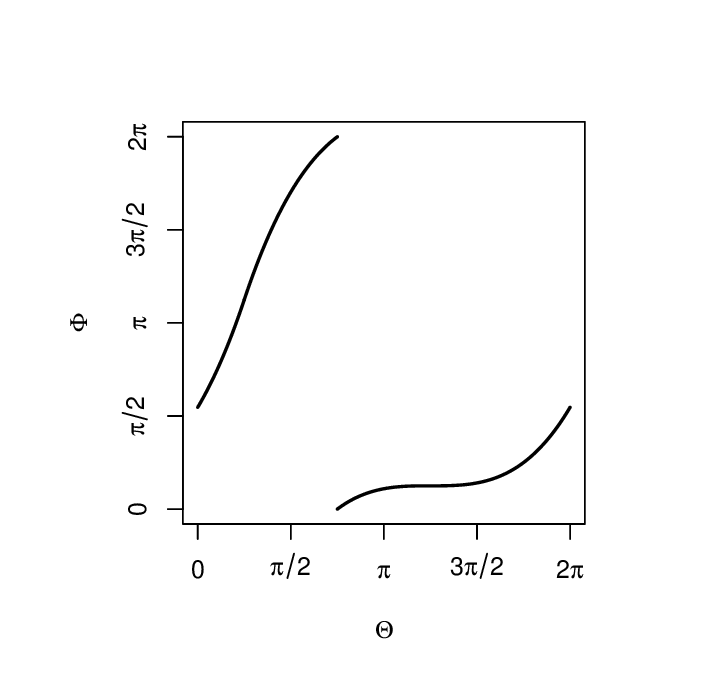} \\
				\includegraphics[width=50mm,height=50mm]{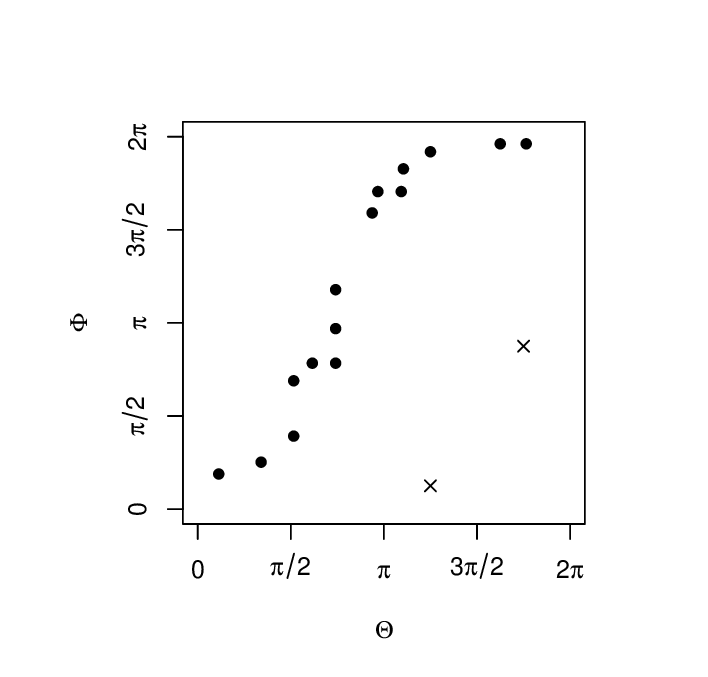}
				\includegraphics[width=50mm,height=50mm]{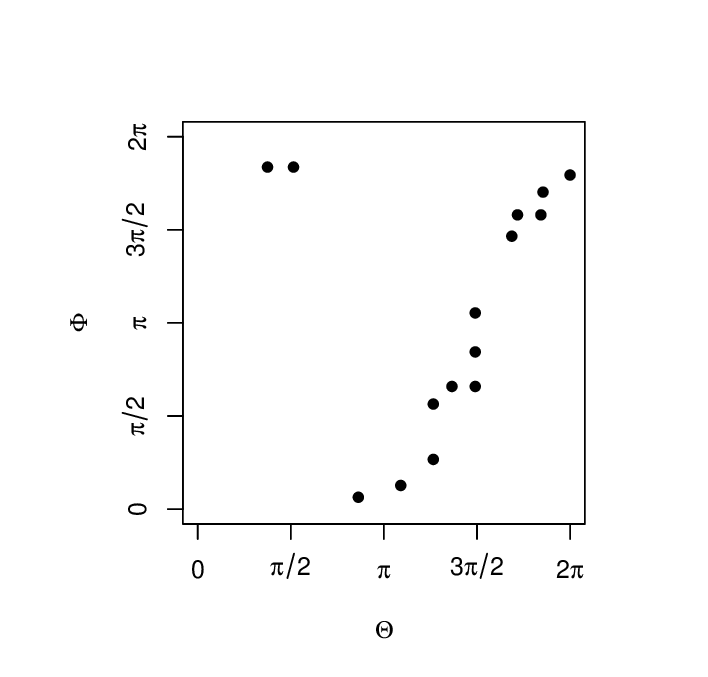}
				\includegraphics[width=50mm,height=50mm]{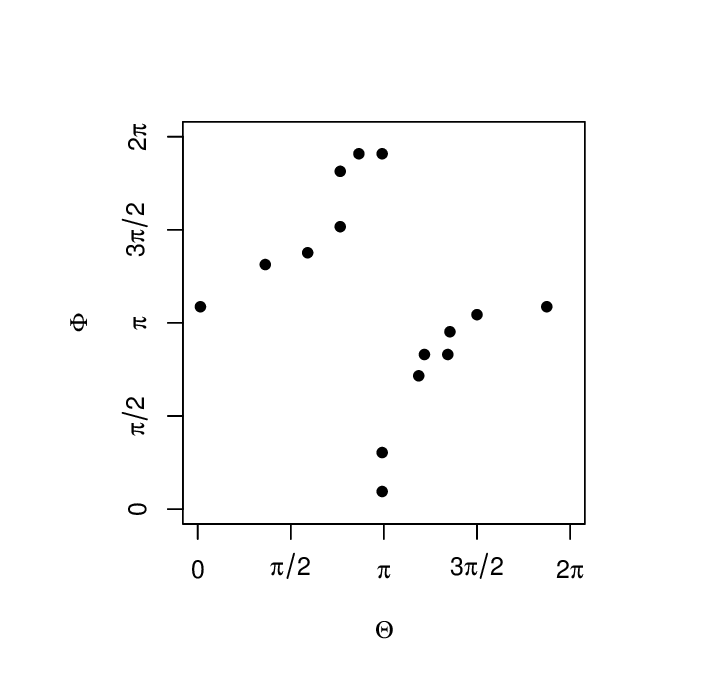} 
		\caption{Supports of perfect positive dependent circular random variables. 
			The upper figures are for a continuous random variable and the lower figures  are for a discrete random variable. The middle and right columns are redrawings of those in the left column when we regard $ (5\pi/4, \pi/8) $ and $ (7\pi/4, 7\pi/8) $ as zero directions (indicated by cross marks in the figures in the left column), respectively.}
		\label{fig:circularnondecreasing}
	\end{figure}

	Now let us consider the pair of circular (angular) random variables $ (\Theta, \Phi)  \in \mathbf{C}^2$. When $ (\Theta, \Phi) $ has perfect positive dependence, its support should be like the examples shown in the left column of Figure \ref{fig:circularnondecreasing}. The upper figures are for a continuous circular random variable and the lower figures are for a discrete variable. The figures in the middle and right columns are redrawings of those in the left column when we regard $ (5\pi/4, \pi/8) $ and $ (7\pi/4, 7\pi/8) $ as zero directions (indicated by cross marks in the figures in the left column), respectively. Due to the arbitrary nature of the zero directions, the all support plots in Figure \ref{fig:circularnondecreasing} should indicate perfect positive dependence.

	Now, we give the equivalence class of the circular Fr\'{e}chet-Hoeffding copula upper bound in the following theorem. 
	\begin{thm}[Circular Fr\'{e}chet-Hoeffding copula upper bound]\label{thm:circular upper bound}
		Let $ 0 \le a \le 1 $. The equivalence class of the circular Fr\'{e}chet-Hoeffding copula upper bound is given by  
		\begin{align*}
		M_a(u,v) 
		= \left\{
		\begin{array}{ll}
		\min(u,v-a) &  (u,v) \in [0,1-a] \times [a,1]  \\
		\min(u+a-1,v) &  (u,v) \in [1-a,1] \times [0,a]  \\
		\max(u+v-1,0) &  \mbox{otherwise} 
		\end{array}
		\right. .
		\end{align*}
	\end{thm} 
	Proof is given in the Appendix. 
	
	The copula $ M_a(u,v) $ is introduced in Exercise 3.9 of \citet{N06}. It is the joint distribution function when the probability mass is uniformly distributed on two line segments, one joining $ (0,a) $ to $ (1-a,1) $ with mass $ 1-a $, and the other joining $ (1-a,0) $ to $ (1,a) $ with mass $ a $. Theorem \ref{thm:circular upper bound} clarifies this $ M_a $ corresponds to the circular Fr\'{e}chet-Hoeffding copula upper bound.  
	
	We can also find the equivalence class of the circular Fr\'{e}chet-Hoeffding copula lower bound in the following theorem. 
	\begin{thm}[Circular Fr\'{e}chet-Hoeffding copula lower bound]\label{thm:circular lower bound}
		Let $ 0 \le a \le 1 $. The equivalence class of the circular Fr\'{e}chet-Hoeffding copula lower bound is given by  
		\begin{align*}
		W_a(u,v) 
		= \left\{
		\begin{array}{ll}
		\max(u+v-a,0) &  (u,v) \in [0,a]^2  \\
		\max(u+v-1,a) &  (u,v) \in [a,1]^2  \\
		\min(u,v) &  \mbox{otherwise} 
		\end{array}
		\right.
		\end{align*}
	\end{thm} 
	Proof is omitted because it is similar to that of Theorem \ref{thm:circular upper bound}. 
	
	The copula $ W_a(u,v) $ is introduced in Example 3.4 of \citet{N06}. It is the joint distribution function when the probability mass is uniformly distributed on two line segments, one joining $ (0,a) $ to $ (a,0) $ with mass $ a $, and the other joining $ (a,1) $ to $ (1,a) $ with mass $ 1-a $. 
	
\begin{figure}[h]
		\centering
			\includegraphics[width=50mm,height=50mm]{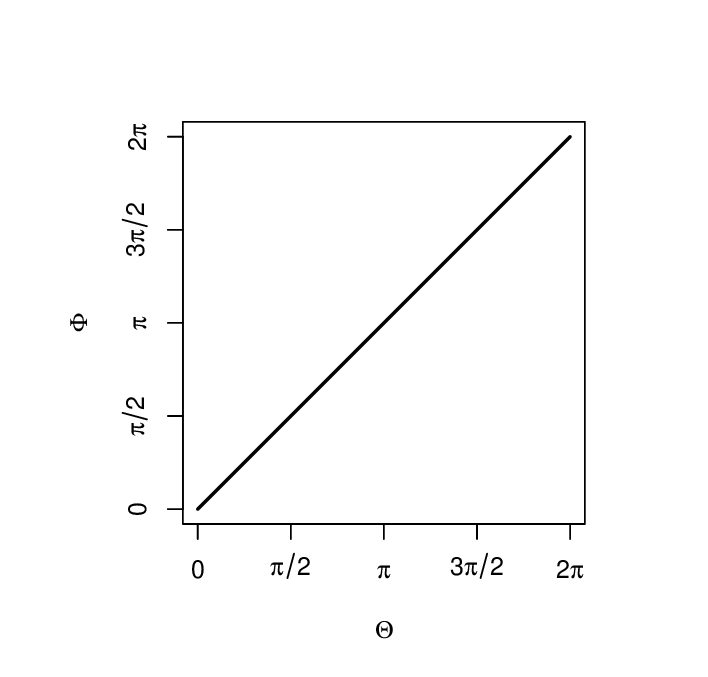} 
			\includegraphics[width=50mm,height=50mm]{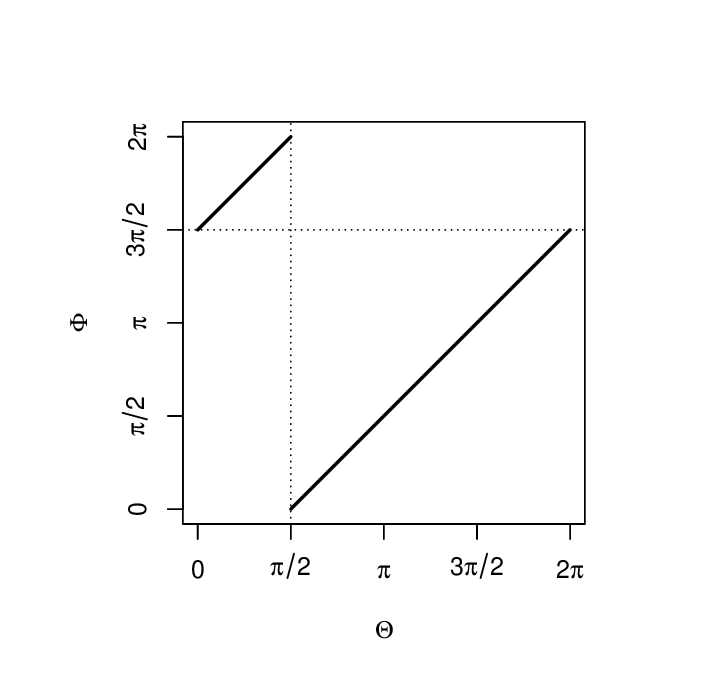} 
			\includegraphics[width=50mm,height=50mm]{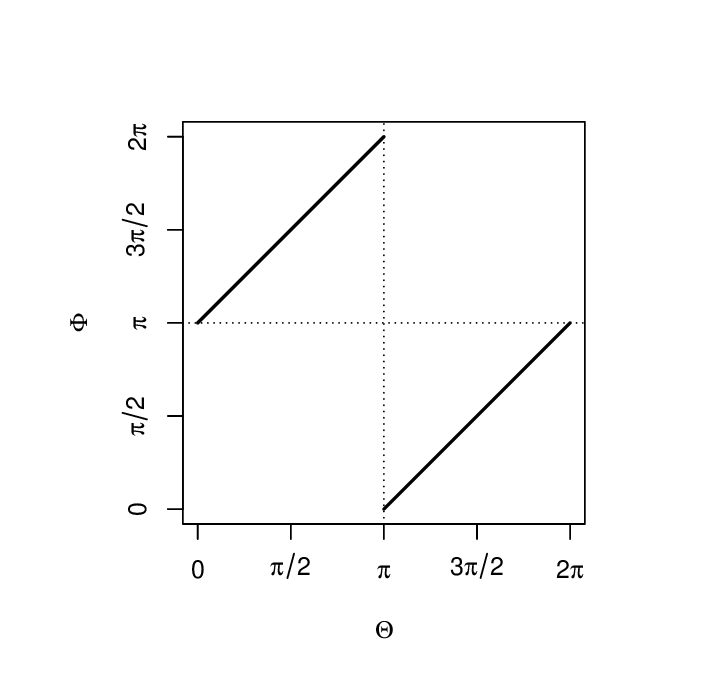} \\
			\includegraphics[width=50mm,height=50mm]{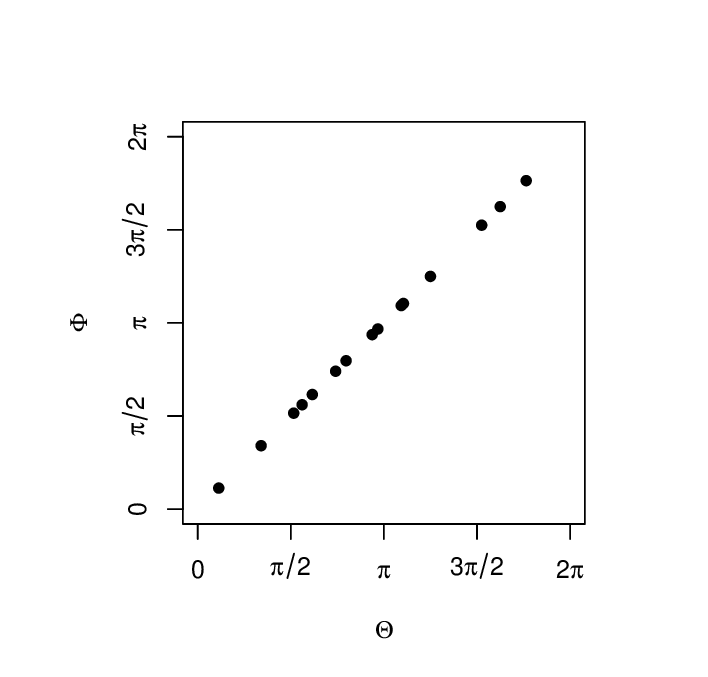} 
			\includegraphics[width=50mm,height=50mm]{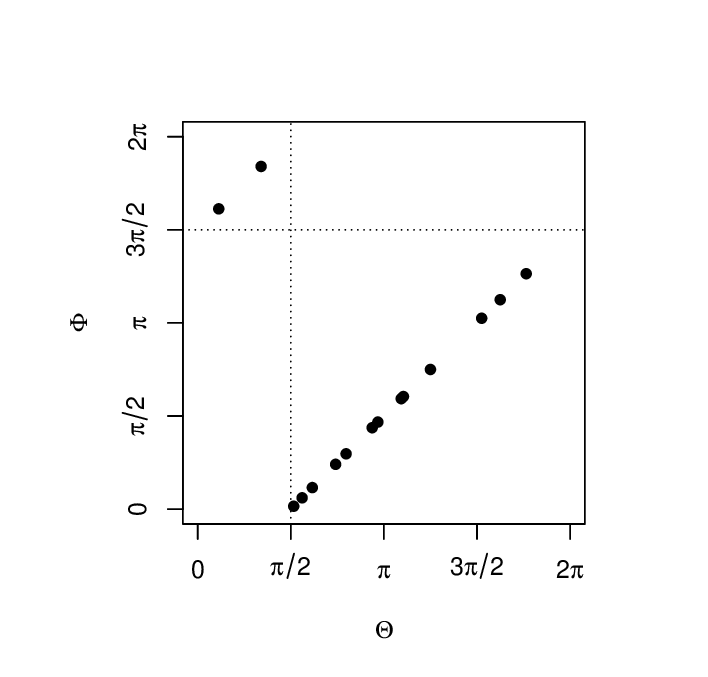} 
			\includegraphics[width=50mm,height=50mm]{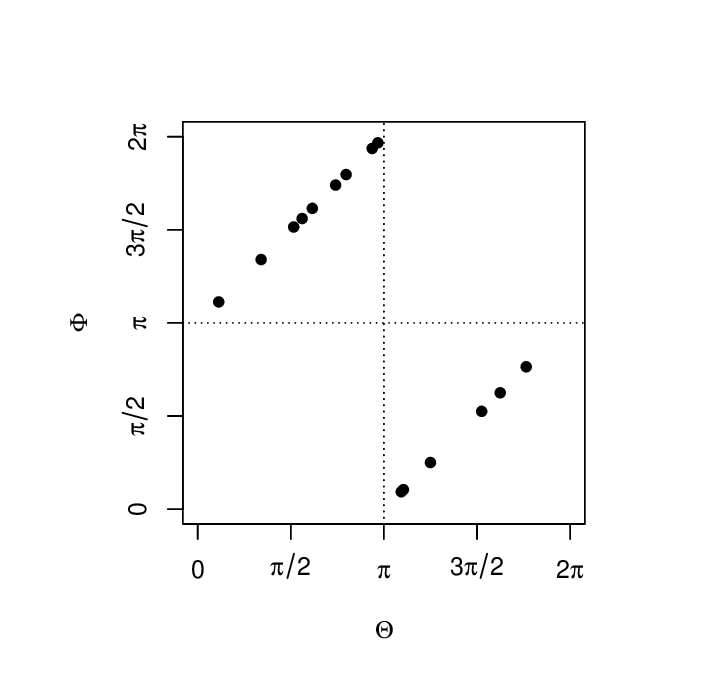} 
		\caption{Perfect positive dependence in the sense of \cite{FL83} with $ \alpha_0=0 $ (left), $ \alpha_0=3\pi/2 $ (middle), $ \alpha_0=\pi $ (right). Upper is for continuous and lower is for discrete.}
		\label{fig:Perfect dependence}
	\end{figure}

	When describing the dependency between $ \Theta $ and $ \Phi $, their periodic structure must be considered. \citet{FL83} proposed the complete dependence between $ \Theta $ and $ \Phi $  as
	\begin{align}
	&\Phi \equiv \Theta + \alpha_0 \qquad \mod (2\pi), \qquad \mbox{positive association} \label{sq:complete positive dependence} \\
	&\Phi \equiv -\Theta + \alpha_0 \qquad \mod (2\pi), \qquad \mbox{negative association} \label{sq:complete negative dependence}
	\end{align}
	where $ \alpha_0 $ is an arbitrary fixed direction. If we use $ \mathbf{C}^2 $ display, the supports of (\ref{sq:complete positive dependence}) 
	are expressed as in Figure \ref{fig:Perfect dependence}. 
	Because of arbitrary nature of $ \alpha_0 $, not only the left column but also the middle and right columns indicate complete positive dependence. Figure \ref{fig:Perfect dependence}  
	is the special case of Figure \ref{fig:circularnondecreasing}. Therefore, the circular perfect dependence with the concept of nondecreasing (nonincreasing) set can be considered as the generalization of complete dependence in the sense of \cite{FL83}.

	\section{Monte Carlo Simulations}\label{sec:Monte Carlo Simulations and Numerical Examples}
	
	In this section, we consider the copula
	\begin{align}\label{eq:Mardia model}
	C_{\gamma,a,b}(u,v)
	= \frac{\gamma^2(1+\gamma)}{2} M_{a}(u,v) + (1-\gamma^2)  \Pi(u,v) + \frac{\gamma^2(1-\gamma)}{2} W_{b}(u,v)& \\ 
	(\gamma,a,b) \in [-1,1] \times \mathbf{I} \times \mathbf{I}.& \nonumber
	\end{align}
	Here, $ M_{a}(u,v) $, $ W_{b}(u,v) $ are the copulas introduced in Theorems \ref{thm:circular upper bound}, \ref{thm:circular lower bound}, and $ \Pi(u,v):=uv $ is an	 independent copula. This is a linear combination of these three copulas and an analog to the model in \citet{M70}. The parameter $ \gamma $ controls the weights of these three copulas, and $ \gamma=1, 0, -1 $ correspond $ M_a $ (perfect positive dependence), $ \Pi $ (independent), $ W_b $ (perfect negative dependence), respectively.  
	
\begin{figure}[h]
		\centering
			\includegraphics[width=50mm,height=50mm]{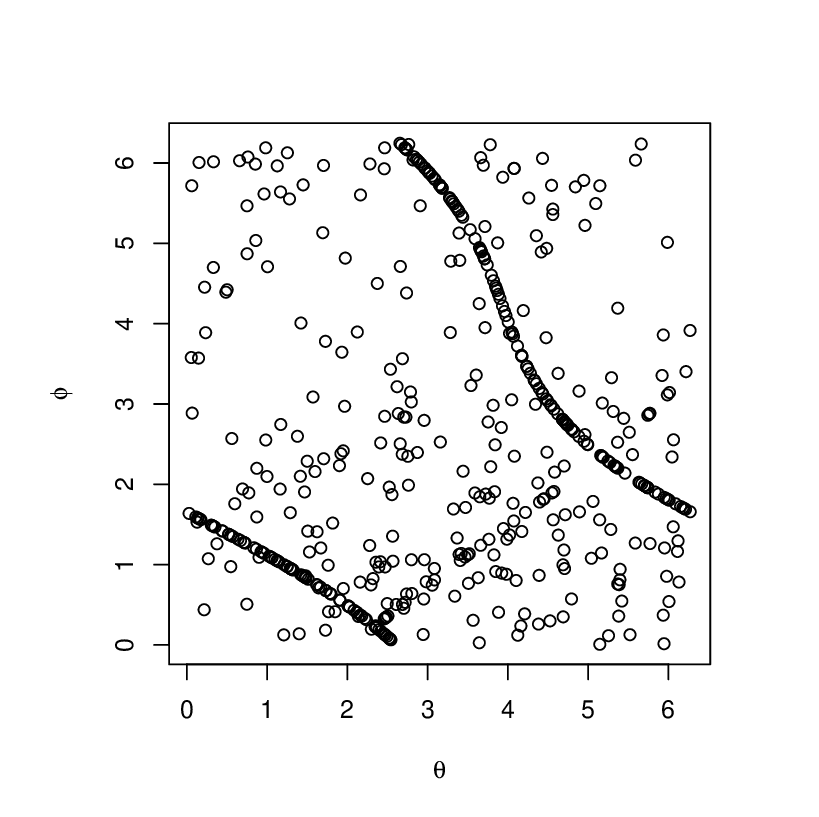}
			\includegraphics[width=50mm,height=50mm]{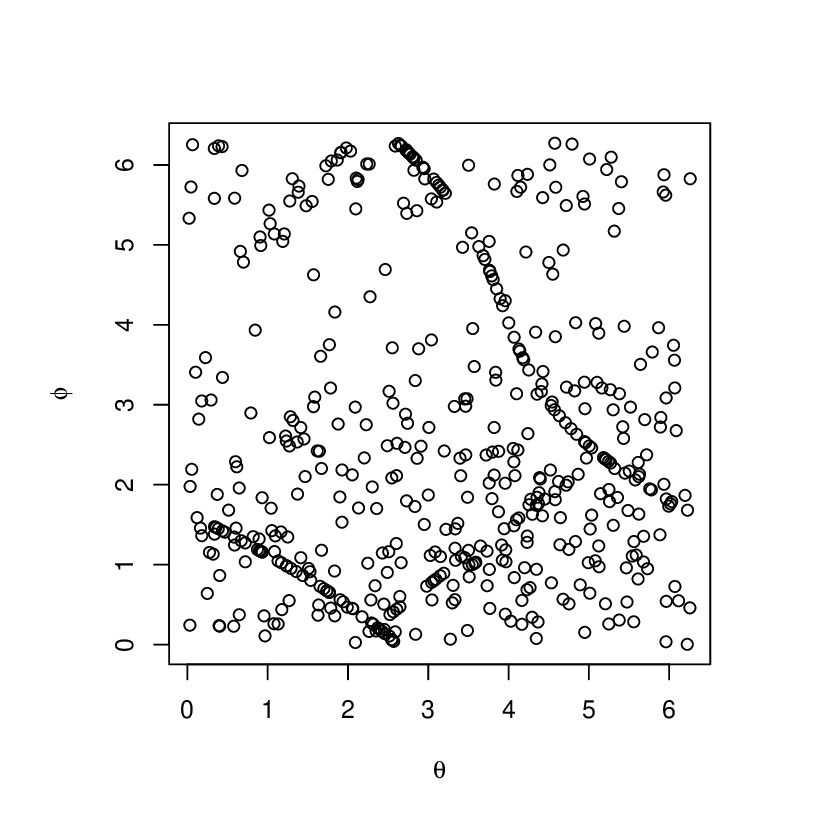}
			\includegraphics[width=50mm,height=50mm]{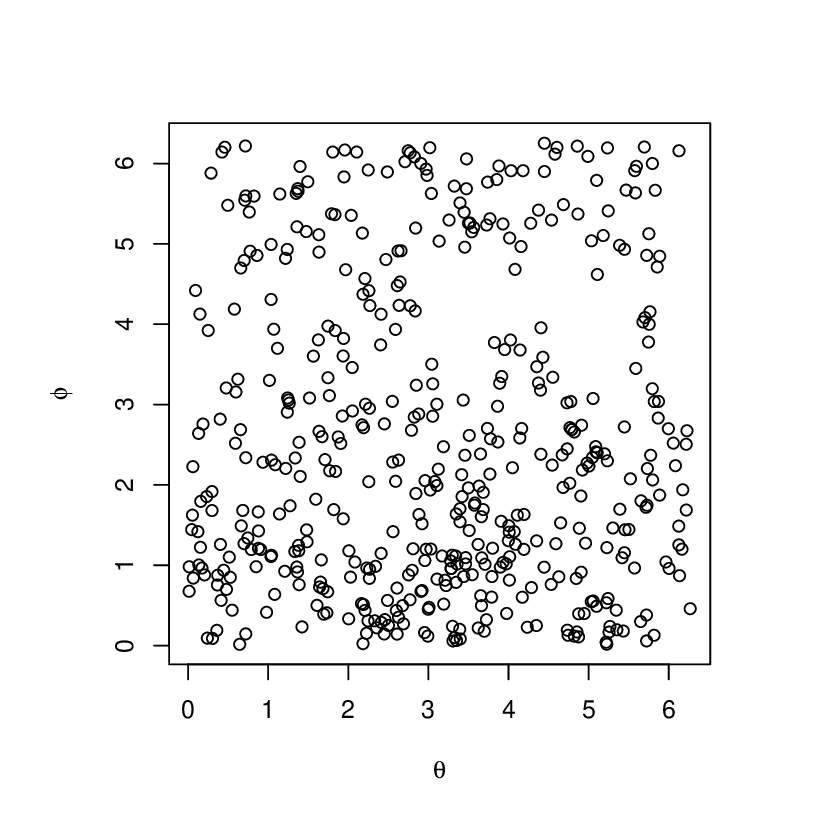} \\
			\includegraphics[width=50mm,height=50mm]{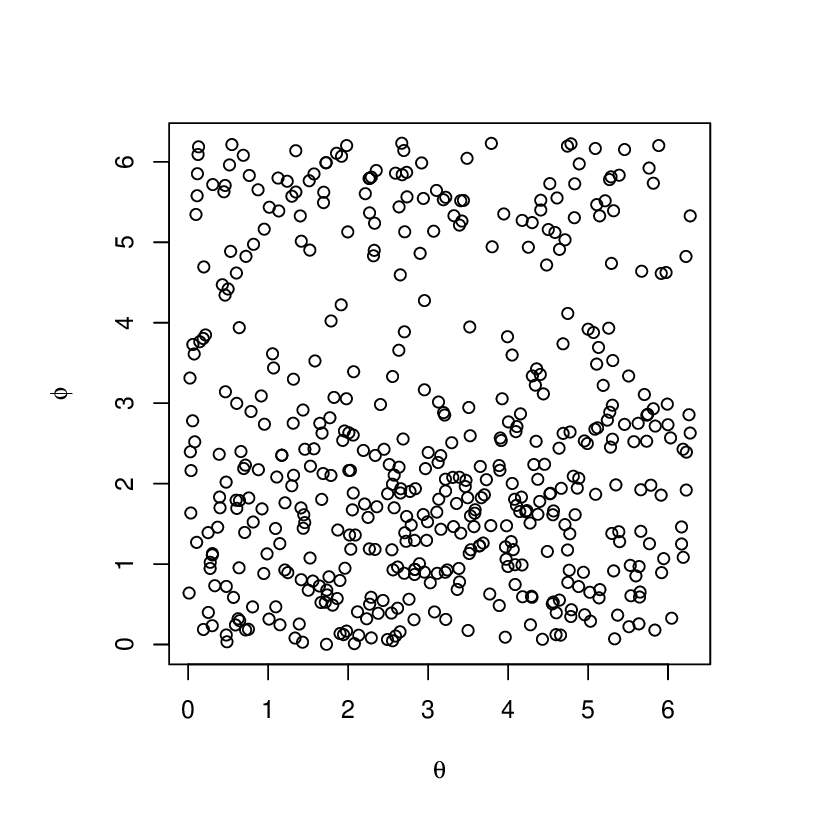}
			\includegraphics[width=50mm,height=50mm]{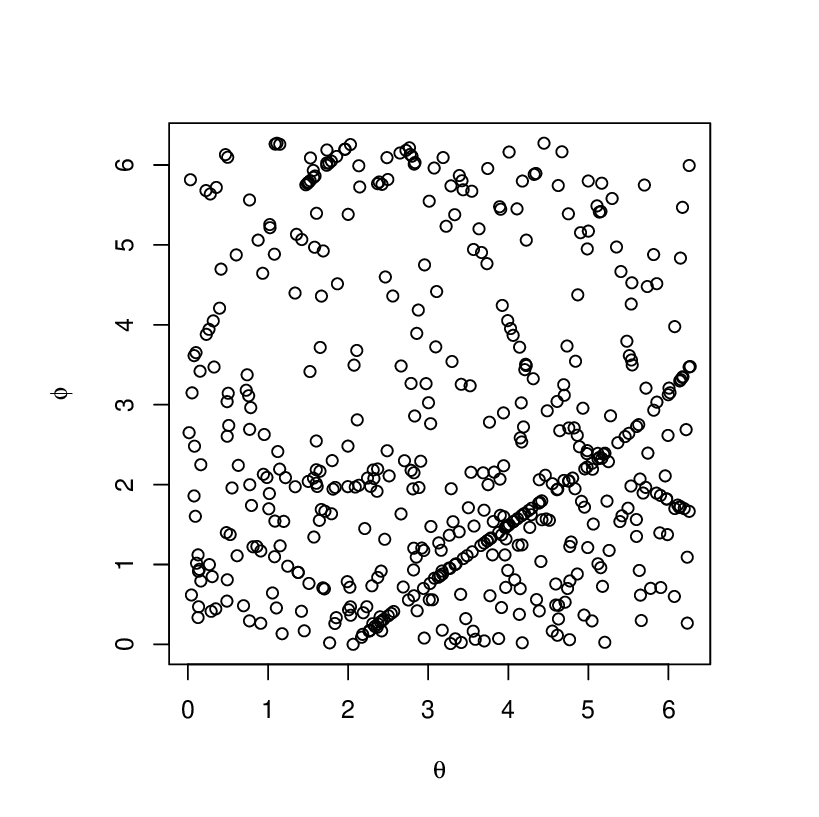}
			\includegraphics[width=50mm,height=50mm]{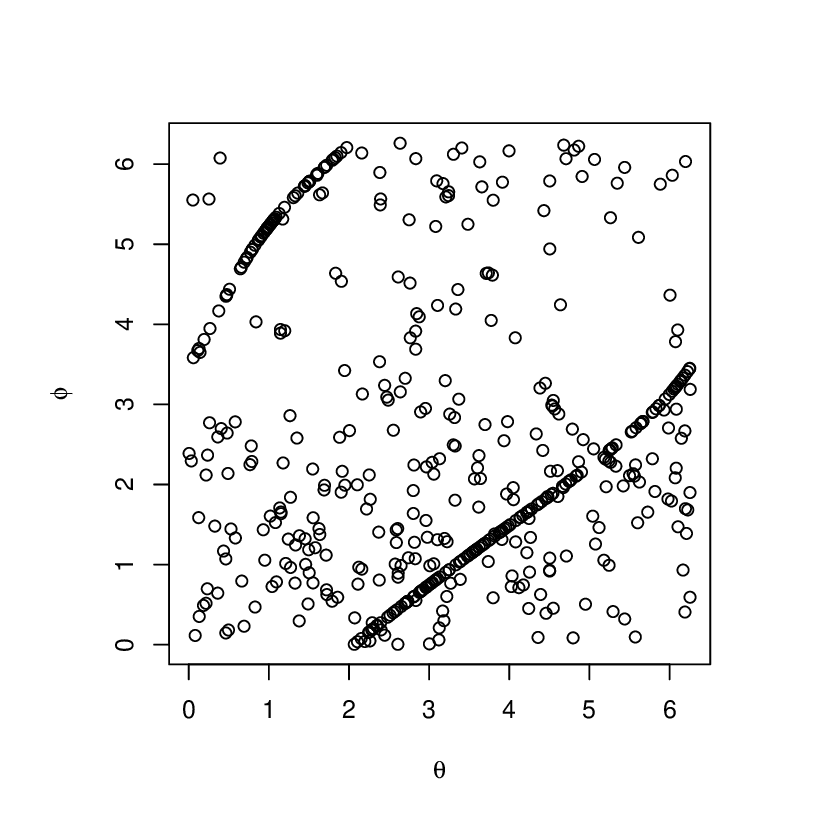}
		\caption{Simulated circular bivariate plots from the model (\ref{eq:Mardia model}). The sample size is $ 500 $. The marginal for $ \theta $ is a Cardioid distribution with $ \rho_F=0.1 $, $ \mu_F=\pi $, and the marginal for $ \phi $ is a Cardioid distribution with $ \rho_G=0.3 $, $ \mu_G=\pi/3 $. The parameters $ a $ and $ b $ are fixed to $ 0.7 $ and $ 0.4 $, respectively. The parameter $ \gamma $ is set to $ -0.7 $ (upper left), $ -0.5 $ (upper middle), $ -0.3 $ (upper right), $ 0.3 $ (lower left), $ 0.5 $ (lower middle), $ 0.7 $ (lower right).}
		\label{fig:simulation}
	\end{figure}

	Now we generate a random circular bivariate sample $ (\theta_i,\phi_i)_{i=1,\ldots,500} $ from (\ref{eq:Mardia model}). Both marginals are set to a Cardioid distribution, whose distribution function is given by
	\begin{align*}
	F_{\mathrm{Ca}}(\theta)=\frac{\rho}{\pi}\sin(\theta-\mu)+\frac{\theta}{2\pi}+\frac{\rho}{\pi}\sin\mu. 
	\end{align*}
	The parameters for marginals are set to $ \rho_F=0.1 $, $ \mu_F=\pi $ for $ \theta $ and $ \rho_G=0.3 $, $ \mu_G=\pi/3 $ for $ \phi $. Figure \ref{fig:simulation} shows the simulated circular bivariate plots with $ \gamma=-0.7 $ (upper left), $ -0.5 $ (upper middle), $ -0.3 $ (upper right), $ 0.3 $ (lower left), $ 0.5 $ (lower middle), $ 0.7 $ (lower right). When $ \gamma $ is close to $ 1 $ and $ -1 $, (\ref{eq:Mardia model}) is close to $ M_a $ and $ W_b $, which are singular copulas. When $ \gamma=\pm 0.7 $, we can clearly see the singular components in the plots. 
	
	\section{Summary and Conclusions}\label{sec:Summary and Conclusions}
	We considered the equivalence class of univariate and bivariate circular distribution functions ascribed to the arbitrary nature of the origin point. Then, we introduced the equivalence class of circular copula functions with Sklar's theorem.    
	Using the concept of the equivalence class of circular copula functions, we introduced circular analogues of Fr\'{e}chet-Hoeffding copula upper and lower bounds. We explained they are the generalizations of complete positive and negative dependence in the sense of \cite{FL83}. 
	We also introduced the circular analogue of Mardia's copula model and simulated a dataset following this model. When the model is close to its extremes, we can find the singular components in the plots.

	\section*{Appendix}
	\paragraph{Proof of Theorem \ref{thm:circular upper bound}}
	Let
	\begin{align*}
	H(\theta,\phi)=M(F(\theta),G(\phi))=\min(F(\theta),G(\phi)). 
	\end{align*}
	Fix arbitrary zero directions $ (\alpha,\beta) \in \mathbf{C}$. From the discussion in Section \ref{sec:Equivalence}, the circular Fr\'{e}chet-Hoeffding copula upper bound is
	\begin{align*}
	&C_{\alpha,\beta}(u,v) \\
	=&H_{\alpha,\beta}(F^{(-1)}_{\alpha}(u),G^{(-1)}_{\beta}(v)) \\
	=&\tilde{H}(F^{(-1)}_{\alpha}(u)+\alpha,G^{(-1)}_{\beta}(v)+\beta)-\tilde{H}(\alpha,G^{(-1)}_{\beta}(v)+\beta)-\tilde{H}(F^{(-1)}_{\alpha}(u)+\alpha,\beta)+\tilde{H}(\alpha,\beta).
	\end{align*}
	When $ (\Theta, \Phi) $ is discrete, $ (u,v) $ is restricted on $ \mathrm{Ran}F_{\alpha} \times \mathrm{Ran}G_{\beta} $.   Now, let us divide $ \mathbf{I}^2 $ into the following four regions:
	\begin{align*}
	&\mathbf{I}_1 = \{ (u,v) \ | \ F^{(-1)}_{\alpha}(u)+\alpha< 2\pi, \ G^{(-1)}_{\beta}(v)+\beta < 2\pi \}
	= \{ (u,v) \ | \ u< 1-F(\alpha), \ v < 1-G(\beta) \} \\
	&\mathbf{I}_2 = \{ (u,v) \ | \ F^{(-1)}_{\alpha}(u)+\alpha< 2\pi, \ G^{(-1)}_{\beta}(v)+\beta \ge 2\pi \} 
	= \{ (u,v) \ | \ u< 1-F(\alpha), \ v \ge 1-G(\beta) \} \\
	&\mathbf{I}_3 = \{ (u,v) \ | \ F^{(-1)}_{\alpha}(u)+\alpha \ge 2\pi, \ G^{(-1)}_{\beta}(v)+\beta < 2\pi \} 
	= \{ (u,v) \ | \ u\ge 1-F(\alpha), \ v < 1-G(\beta) \}\\
	&\mathbf{I}_4 = \{ (u,v) \ | \ F^{(-1)}_{\alpha}(u)+\alpha \ge 2\pi, \ G^{(-1)}_{\beta}(v)+\beta \ge 2\pi \} 
	= \{ (u,v) \ | \ u\ge 1-F(\alpha), \ v \ge 1-G(\beta) \}
	\end{align*}
	
	First, for $ (u,v) \in \mathbf{I}_1 $, 
	\begin{align*}
	&C_{\alpha,\beta}(u,v) \\
	=&H(F^{(-1)}_{\alpha}(u)+\alpha,G^{(-1)}_{\beta}(v)+\beta)-H(\alpha,G^{(-1)}_{\beta}(v)+\beta)-H(F^{(-1)}_{\alpha}(u)+\alpha,\beta)+H(\alpha,\beta) \\
	=&\min\{F(F^{(-1)}_{\alpha}(u)+\alpha),G(G^{(-1)}_{\beta}(v)+\beta)\}
	-\min\{F(\alpha),G(G^{(-1)}_{\beta}(v)+\beta)\}\\
	&-\min\{F(F^{(-1)}_{\alpha}(u)+\alpha),G(\beta)\}
	+\min\{F(\alpha),G(\beta)\}
	\end{align*}
	Here, 
	\begin{align*}
	&F(F^{(-1)}_{\alpha}(u)+\alpha)
	=F_\alpha(F^{(-1)}_{\alpha}(u))+F(\alpha)
	=u+F(\alpha) \\
	&G(G^{(-1)}_{\beta}(v)+\beta)
	=G_\beta(G^{(-1)}_{\beta}(v))+G(\beta)
	=v+G(\beta).
	\end{align*}
	Therefore,
	\begin{align*}
	C_{\alpha,\beta}(u,v) 
	=&\min\{u+F(\alpha),v+G(\beta)\}
	-\min\{F(\alpha),v+G(\beta)\} \\
	&-\min\{u+F(\alpha),G(\beta)\}
	+\min\{F(\alpha),G(\beta)\}.
	\end{align*}
	When $ G(\beta) > F(\alpha) $, we divide $ \mathbf{I}_1 $ into 
	\begin{align*}
	&\mathbf{I}^G_{1-1} = \{ (u,v) \in \mathbf{I}_1 \ | \ u \le G(\beta)-F(\alpha) \} \\
	&\mathbf{I}^G_{1-2} = \{ (u,v) \in \mathbf{I}_1 \ | \ u > G(\beta)-F(\alpha), \ v > u-(G(\beta)-F(\alpha)) \} \\
	&\mathbf{I}^G_{1-3} = \{ (u,v) \in \mathbf{I}_1 \ | \ v \le u-(G(\beta)-F(\alpha)) \}.
	\end{align*}
	When $ G(\beta) \le F(\alpha) $, we divide $ \mathbf{I}_1 $ into 
	\begin{align*}
	&\mathbf{I}^F_{1-1} = \{ (u,v) \in \mathbf{I}_1 \ | \ v \le F(\alpha)-G(\beta) \} \\
	&\mathbf{I}^F_{1-2} = \{ (u,v) \in \mathbf{I}_1 \ | \ v > F(\alpha)-G(\beta), \ v \le u+(F(\alpha)-G(\beta)) \} \\
	&\mathbf{I}^F_{1-3} = \{ (u,v) \in \mathbf{I}_1 \ | \ v > u+(F(\alpha)-G(\beta)) \}. 
	\end{align*}
	Then, 
	\begin{align}
	\label{eq:G1-1}
	\begin{split}
	C_{\alpha,\beta}(u,v) 
	= \left\{ 
	\begin{array}{cl}
	0 & \quad (u,v) \in \mathbf{I}^G_{1-1} \\[1mm]
	u-(G(\beta)-F(\alpha)) & \quad (u,v) \in \mathbf{I}^G_{1-2} \\[1mm]
	v & \quad (u,v) \in \mathbf{I}^G_{1-3} \\[1mm]
	0 & \quad (u,v) \in \mathbf{I}^F_{1-1} \\[1mm]
	v-(F(\alpha)-G(\beta)) & \quad (u,v) \in \mathbf{I}^F_{1-2} \\[1mm]
	u & \quad (u,v) \in \mathbf{I}^F_{1-3} \\
	\end{array}
	\right.    
	\end{split}
	\end{align}
	
	Second, for $ (u,v) \in \mathbf{I}_2 $, 
	\begin{align*}
	&C_{\alpha,\beta}(u,v) \\
	=&[F(\alpha+F^{(-1)}_{\alpha}(u))+H(F^{(-1)}_{\alpha}(u)+\alpha,G^{(-1)}_{\beta}(v)+\beta-2\pi)] \\
	&-[F(\alpha)+H(\alpha, G^{(-1)}_{\beta}(v)+\beta-2\pi)]
	-H(F^{(-1)}_{\alpha}(u)+\alpha,\beta)
	+H(\alpha,\beta) \\
	=&[F(\alpha+F^{(-1)}_{\alpha}(u))+\min\{F(F^{(-1)}_{\alpha}(u)+\alpha), G(G^{(-1)}_{\beta}(v)+\beta-2\pi) \}] \\
	&-[F(\alpha)+\min\{F(\alpha), G(G^{(-1)}_{\beta}(v)+\beta-2\pi)\}]
	-\min\{F(F^{(-1)}_{\alpha}(u)+\alpha), G(\beta) \} \\
	&+\min\{F(\alpha), G(\beta) \} 
	\end{align*}
	Here, 
	\begin{align*}
	&F(F^{(-1)}_{\alpha}(u)+\alpha)
	=F_\alpha(F^{(-1)}_{\alpha}(u))+F(\alpha)
	=u+F(\alpha) \\
	&G(G^{(-1)}_{\beta}(v)+\beta-2\pi)
	=G_\beta(G^{(-1)}_{\beta}(v))-G_\beta(2\pi-\beta)
	=v-(1-G(\beta)). 
	\end{align*}
	Therefore,
	\begin{align*}
	C_{\alpha,\beta}(u,v) 
	=&u+\min\{u+F(\alpha),v-(1-G(\beta))\}
	-\min\{F(\alpha),v-(1-G(\beta))\} \\
	&-\min\{u+F(\alpha),G(\beta)\}
	+\min\{F(\alpha),G(\beta)\}
	\end{align*}
	When $ G(\beta) > F(\alpha) $, we divide $ \mathbf{I}_2 $ into 
	\begin{align*}
&\mathbf{I}^G_{2-1} = \{ (u,v) \in \mathbf{I}_2 \ | \ u \le G(\beta)-F(\alpha), \ v \le 1-(G(\beta)-F(\alpha)) \} \\
&\mathbf{I}^G_{2-2} = \{ (u,v) \in \mathbf{I}_2 \ | \ u > G(\beta)-F(\alpha), \ v \le 1-(G(\beta)-F(\alpha)) \} \\
&\mathbf{I}^G_{2-3} = \{ (u,v) \in \mathbf{I}_2 \ | \ v > u+\{1-(G(\beta)-F(\alpha)) \} \, \} \\
&\mathbf{I}^G_{2-4} = \{ (u,v) \in \mathbf{I}_2 \ | \ u \le G(\beta)-F(\alpha), \  1-(G(\beta)-F(\alpha)) < v \le u+\{1-(G(\beta)-F(\alpha)) \} \,  \} \\
&\mathbf{I}^G_{2-5} = \{ (u,v) \in \mathbf{I}_2 \ | \ u > G(\beta)-F(\alpha), \ v > 1-(G(\beta)-F(\alpha)) \}. 
\end{align*}
	When $ G(\beta) \le F(\alpha) $, we do not divide $ \mathbf{I}_2 $ further and rename it as $ \mathbf{I}^F_{2-1} $. Then, 
	\begin{align}
	\label{eq:G2-1}
	\begin{split}
	C_{\alpha,\beta}(u,v) 
	= \left\{ 
	\begin{array}{cl}
	0 & \quad (u,v) \in \mathbf{I}^G_{2-1} \\[1mm]
	u-(G(\beta)-F(\alpha)) & \quad (u,v) \in \mathbf{I}^G_{2-2} \\[1mm]
	u & \quad (u,v) \in \mathbf{I}^G_{2-3} \\[1mm]
	v-\{1-(G(\beta)-F(\alpha))\} & \quad (u,v) \in \mathbf{I}^G_{2-4} \\[1mm]
	u+v-1 & \quad (u,v) \in \mathbf{I}^G_{2-5} \\[1mm]
	u & \quad (u,v) \in \mathbf{I}^F_{2-1} \\
	\end{array}
	\right.    
	\end{split}
	\end{align}
	
	Third, for $ (u,v) \in \mathbf{I}_3 $, 
	\begin{align*}
	&C_{\alpha,\beta}(u,v) \\
	=&[G(\beta+G^{(-1)}_{\beta}(v))+H(F^{(-1)}_{\alpha}(u)+\alpha-2\pi, G^{(-1)}_{\beta}(v)+\beta)] \\
	&-[G(\beta)+H(F^{(-1)}_{\alpha}(u)+\alpha-2\pi, \beta)]
	-H(\alpha,G^{(-1)}_{\beta}(v)+\beta)
	+H(\alpha,\beta) \\
	=&[G(\beta+G^{(-1)}_{\beta}(v))+\min\{F(F^{(-1)}_{\alpha}(u)+\alpha-2\pi), G(G^{(-1)}_{\beta}(v)+\beta)\}] \\
	&-[G(\beta)+\min\{F(F^{(-1)}_{\alpha}(u)+\alpha-2\pi), G(\beta)\}]
	-\min\{F(\alpha),G(G^{(-1)}_{\beta}(v)+\beta)\} \\
	&+\min\{F(\alpha),G(\beta)\}.
	\end{align*}
	Here, 
	\begin{align*}
	&F(F^{(-1)}_{\alpha}(u)+\alpha-2\pi)
	=F_\alpha(F^{(-1)}_{\alpha}(u))-F_\alpha(2\pi-\alpha)
	=u-(1-F(\alpha)) \\
	&G(G^{(-1)}_{\beta}(v)+\beta)
	=G_\beta(G^{(-1)}_{\beta}(v))+G(\beta)
	=v+G(\beta).
	\end{align*}
	Therefore,
	\begin{align*}
	C_{\alpha,\beta}(u,v) 
	=&v+\min\{u-(1-F(\alpha)),v+G(\beta)\}
	-\min\{u-(1-F(\alpha)),G(\beta)\} \\
	&-\min\{F(\alpha),v+G(\beta)\}
	+\min\{F(\alpha),G(\beta)\}
	\end{align*}
	When $ G(\beta) > F(\alpha) $, we do not divide $ \mathbf{I}_3 $ further and rename it $ \mathbf{I}^G_{3-1} $. 
	When $ G(\beta) \le F(\alpha) $, we divide $ \mathbf{I}_3 $ into 
	\begin{align*}
&\mathbf{I}^F_{3-1} = \{ (u,v) \in \mathbf{I}_3 \ | \ u \le 1-(F(\alpha)-G(\beta)), \ v \le F(\alpha)-G(\beta) \} \\
&\mathbf{I}^F_{3-2} = \{ (u,v) \in \mathbf{I}_3 \ | \ u > 1-(F(\alpha)-G(\beta)), \ v \le F(\alpha)-G(\beta), \ v>u-\{1-(F(\alpha)-G(\beta)) \} \, \} \\
&\mathbf{I}^F_{3-3} = \{ (u,v) \in \mathbf{I}_3 \ | \ v \le u-\{1-(F(\alpha)-G(\beta)) \} \, \} \\
&\mathbf{I}^F_{3-4} = \{ (u,v) \in \mathbf{I}_3 \ | \ u \le 1-(F(\alpha)-G(\beta)), \ v > F(\alpha)-G(\beta)  \} \\
&\mathbf{I}^F_{3-5} = \{ (u,v) \in \mathbf{I}_3 \ | \ u > 1-(F(\alpha)-G(\beta)), \ v > F(\alpha)-G(\beta) \}. 
\end{align*}
Then, 
\begin{align}
\label{eq:G3-1}
\begin{split}
C_{\alpha,\beta}(u,v) 
= \left\{ 
\begin{array}{cl}
v & \quad (u,v) \in \mathbf{I}^G_{3-1} \\[1mm]
0 & \quad (u,v) \in \mathbf{I}^F_{3-1} \\[1mm]
u-\{1-(F(\alpha)-G(\beta))\} & \quad (u,v) \in \mathbf{I}^F_{3-2} \\[1mm]
v & \quad (u,v) \in \mathbf{I}^F_{3-3} \\[1mm]
v-(F(\alpha)-G(\beta)) & \quad (u,v) \in \mathbf{I}^F_{3-4} \\[1mm]
u+v-1 & \quad (u,v) \in \mathbf{I}^F_{3-5} \\
\end{array}
\right.    
\end{split}
\end{align}

	Lastly, for $ (u,v) \in \mathbf{I}_4 $, 
	\begin{align*}
	&C_{\alpha,\beta}(u,v) \\
	=&[1+F(F^{(-1)}_{\alpha}(u)+\alpha-2\pi)+G(G^{(-1)}_{\beta}(v)+\beta-2\pi)+H(F^{(-1)}_{\alpha}(u)+\alpha-2\pi,G^{(-1)}_{\beta}(v)+\beta-2\pi)] \\
	&-[F(\alpha)+H(\alpha,G^{(-1)}_{\beta}(v)+\beta-2\pi)]
	-[G(\beta)+H(F^{(-1)}_{\alpha}(u)+\alpha-2\pi,\beta)]+H(\alpha,\beta) \\
	=&[1+F(F^{(-1)}_{\alpha}(u)+\alpha-2\pi)+G(G^{(-1)}_{\beta}(v)+\beta-2\pi)+\min\{F(F^{(-1)}_{\alpha}(u)+\alpha-2\pi),G(G^{(-1)}_{\beta}(v)+\beta-2\pi)\} ]\\
	&-[F(\alpha)+\min\{F(\alpha),G(G^{(-1)}_{\beta}(v)+\beta-2\pi)\}]
	-[G(\beta)+\min\{F(F^{(-1)}_{\alpha}(u)+\alpha-2\pi),G(\beta)\}]\\
	&+\min\{F(\alpha),G(\beta)\}. 
	\end{align*}
	Here, 
	\begin{align*}
	&F(F^{(-1)}_{\alpha}(u)+\alpha-2\pi)
	=F_\alpha(F^{(-1)}_{\alpha}(u))-F_\alpha(2\pi-\alpha)
	=u-(1-F(\alpha)) \\
	&G(G^{(-1)}_{\beta}(v)+\beta-2\pi)
	=G_\beta(G^{(-1)}_{\beta}(v))-G_\beta(2\pi-\beta)
	=v-(1-G(\beta)). 
	\end{align*}
	Therefore,
	\begin{align*}
	&C_{\alpha,\beta}(u,v) \\
	=&u+v-1+\min\{u-(1-F(\alpha)),v-(1-G(\beta))\} \\
	&-\min\{F(\alpha),v-(1-G(\beta))\}-\min\{u-(1-F(\alpha)),G(\beta)\}
	+\min\{F(\alpha),G(\beta)\}.
	\end{align*}
	When $ G(\beta) > F(\alpha) $, we divide $ \mathbf{I}_4 $ into 
	\begin{align*}
	&\mathbf{I}^G_{4-1} = \{ (u,v) \in \mathbf{I}_4 \ | \ v \le u-(G(\beta)-F(\alpha)) \} \\
	&\mathbf{I}^G_{4-2} = \{ (u,v) \in \mathbf{I}_4 \ | \ v \le 1-(G(\beta)-F(\alpha)), \ v > u-(G(\beta)-F(\alpha)) \} \\
	&\mathbf{I}^G_{4-3} = \{ (u,v) \in \mathbf{I}_4 \ | \ v > 1-(G(\beta)-F(\alpha)) \} 
	\end{align*}
	When $ G(\beta) \le F(\alpha) $, we divide $ \mathbf{I}_4 $ into 
	\begin{align*}
	&\mathbf{I}^F_{4-1} = \{ (u,v) \in \mathbf{I}_4 \ | \ v > u+(F(\alpha)-G(\beta)) \} \\
	&\mathbf{I}^F_{4-2} = \{ (u,v) \in \mathbf{I}_4 \ | \ u \le 1-(F(\alpha)-G(\beta)), \ v \le u+(F(\alpha)-G(\beta)) \} \\
	&\mathbf{I}^F_{4-3} = \{ (u,v) \in \mathbf{I}_4 \ | \ u > 1-(F(\alpha)-G(\beta)) \}. 
	\end{align*}
		Then, 
	\begin{align}
	\label{eq:F4-3}
	\begin{split}
	C_{\alpha,\beta}(u,v) 
	= \left\{ 
	\begin{array}{cl}
	v & \quad (u,v) \in \mathbf{I}^G_{4-1} \\[1mm]
	u-(G(\beta)-F(\alpha)) & \quad (u,v) \in \mathbf{I}^G_{4-2} \\[1mm]
	u+v-1 & \quad (u,v) \in \mathbf{I}^G_{4-3} \\[1mm]
	u & \quad (u,v) \in \mathbf{I}^F_{4-1} \\[1mm]
	v-(F(\alpha)-G(\beta)) & \quad (u,v) \in \mathbf{I}^F_{4-2} \\[1mm]
	u+v-1 & \quad (u,v) \in \mathbf{I}^F_{4-3} \\
	\end{array}
	\right.    
	\end{split}
	\end{align}
	
	From (\ref{eq:G1-1})--(\ref{eq:F4-3}), we find $ C_{\alpha,\beta}(u,v) $ is equivalent to $ M_a(u,v) $ in Theorem \ref{thm:circular upper bound} by taking
	\begin{align*}
	a=\left\{
	\begin{array}{ll}
	1-(G(\beta)-F(\alpha)) & (F(\alpha) \le G(\beta)) \\
	F(\alpha)-G(\beta) & (F(\alpha) > G(\beta)) 
	\end{array}
	\right. .
	\end{align*}

\newpage	
	\section*{Acknowledgements}
	Financial support for this research was received from JSPS KAKENHI in the form of grant 18K11193. 
	We would like to thank Editage (www.editage.jp) for English language editing.

\end{document}